\documentclass[hidelinks]{article}
\usepackage[utf8]{inputenc}
\usepackage{titlesec}
\usepackage[includehead,a4paper,headheight=13.6pt]{geometry} 
\usepackage[indent=15pt]{parskip} 
\usepackage{amssymb,amsmath,amsthm}
\usepackage{indentfirst}
\usepackage{hyperref}
\usepackage[nameinlink]{cleveref}
\usepackage{xcolor}

\newcommand\nnfootnote[1]{%
  \begin{NoHyper}
  \renewcommand\thefootnote{}\footnote{#1}%
  \addtocounter{footnote}{-1}%
  \end{NoHyper}
}

\newcommand{\address}{{
\bigskip\footnotesize
\noindent Iason Efraimidis, \textsc{Departamento de Matem\'aticas, Universidad Aut\'onoma De Madrid, 28049 Madrid, Spain}
\noindent\textit{E-mail address}: \texttt{iason.efraimidis@uam.es}

\noindent Adri\'an Llinares, \textsc{Departamento de Matem\'aticas, Universidad Aut\'onoma De Madrid, 28049 Madrid, Spain}
\noindent\textit{E-mail address}: \texttt{adrian.llinares@uam.es}\\
\noindent\textit{URL:} \texttt{https://sites.google.com/view/adrian-llinares/home}

\noindent Sullivan F.\ MacDonald \textsc{Department of Mathematics, University of Toronto, ON Canada M5S2E4}\\
\noindent\textit{E-mail address}: \texttt{sullivan@math.toronto.edu}

\noindent Ignacio Uriarte-Tuero, \textsc{Department of Mathematics, University of Toronto, ON Canada M5S2E4}\\
\noindent\textit{E-mail address}: \texttt{iuriarte@math.toronto.edu}

\noindent Dragan Vukoti\'c, \textsc{Departamento de Matem\'aticas, Universidad Aut\'onoma De Madrid, 28049 Madrid, Spain}
\noindent\textit{E-mail address}: \texttt{dragan.vukotic@uam.es}\\
\noindent\textit{URL:} \texttt{http://www.uam.es/dragan.vukotic}
}}

\newtheorem{thm}{Theorem}
\newtheorem{lem}[thm]{Lemma}
\newtheorem{cor}{Corollary}[thm]

\newcommand{\re}{\mathrm{Re}\,}
\newcommand{\im}{\mathrm{Im}\,}
\newcommand{\rmk}{\noindent\textit{Remark}: }

\titleformat*{\section}{\large\bfseries}
\titleformat*{\subsection}{\normalsize\bfseries}

\title{A property of boundary extensions of atomic singular inner functions}

\author{\normalsize Iason Efraimidis, Adri\'an Llinares, Sullivan MacDonald, Ignacio Uriarte-Tuero, Dragan Vukoti\'c}
\date{}

\begin{document}
\maketitle
\nnfootnote{\textit{Date}: November 4, 2026}
\nnfootnote{\textit{2020 Mathematics Subject Classification}. Primary 30J15, Secondary 28A20}
\nnfootnote{\textit{Key words and phrases}. Singular inner function, distribution functions, measure preserving transformations}

\begin{abstract}
Using a property of the Hilbert transform, we prove a measure theoretic result concerning the extensions of atomic singular inner functions to the circle. Various alternative proofs are presented using the mean value theorem for harmonic functions and a well-known property of Blaschke products. Using the main result, we establish several facts concerning atomic singular inner functions and their logarithms. 
\end{abstract}

\section{Introduction and main results}

Let \(\mathrm{Hol}(\mathbb D)\) denote the set of holomorphic functions on the unit disc in \(\mathbb C\). An atomic singular inner function is a function \(f\in\mathrm{Hol}(\mathbb D)\) of the form
\begin{equation}\label{ASIF}
    f(z)=\eta \exp\bigg(-\sum_{j=1}^N\lambda_j\frac{1+\alpha_jz}{1-\alpha_jz}\bigg),
\end{equation}
where \(\eta\in\partial\mathbb D\), \(N\in\mathbb N\), \(\alpha_1,\dots,\alpha_N\) are distinct points on \(\partial\mathbb D\), and \(\lambda_1,\dots,\lambda_N>0\). Equivalently,
\[
    f(z)=\eta \exp\bigg(-\frac{1}{2\pi }\int_\mathbb TH(z,\theta)d\mu(\theta)\bigg), \qquad H(z,\theta)=\frac{e^{i\theta}+z}{e^{i\theta}-z},
\]
where \(H\) is the Herglotz kernel and \(\mu\) is a discrete atomic measure supported at finitely many points \(\theta_1,\dots,\theta_N\in  \mathbb T\),
\[
    \mu=\sum_{j=1}^N\lambda_j\delta_{\theta_j}.
\]
The points \(\alpha_j\)  and \(\theta_j\) are related by the equations \(\alpha_j=e^{-i\theta_j}\), up to reindexing. 

Functions of the form \eqref{ASIF} appear in various settings. They arise as extremal functions for non-linear coefficient problems in complex analysis, such as the Krzy\.z conjecture \cite{HummelScheinbergZalcman1977}, and as factors of extremals for the Hummel-Scheinberg-Zalcman conjecture and related problems \cite{Suffridge1990}. Recently, Instanes, Llinares, \& Vukoti\'c \cite{InstanesLlinaresVukotic2026} showed that extremal functions for a much more general class of coefficient problems are of the form \eqref{ASIF}. 

While these applications are our primary motivation for studying atomic singular inner functions, we note that singular inner functions appear also in the canonical factorization theorem for Nevanlinna functions (see e.g.\ \cite[Thm.\ 5.5]{Garnett2007}), which states that if \(f\) belongs to the Nevanlinna class (i.e. it satisfies \(\sup_{0<r<1}\int_\mathbb T\log^+ |f(re^{i\theta})|d\theta<\infty\)), then it can be written \(f=\eta BFS_1/S_2\), where \(\eta\in\partial\mathbb D\), \(B\) is a Blaschke product, \(F\) is an outer function, and \(S_1\) and \(S_2\) are singular inner functions. This generalizes the well-known Beurling factorization for \(H^p\) functions. Our results may therefore be useful in various settings where these factorizations appear. 

Fixing \(f\) as in \eqref{ASIF}, we observe that \(f\in H^\infty(\mathbb D)\), and indeed \(\sup_\mathbb D|f|\leq 1\) since \(\re \log f<0\) everywhere in \(\mathbb D\). In what follows, we take \(\eta=1\). As \(f\) is non-vanishing in \(\mathbb D\), it has a holomorphic logarithm in \(\mathbb D\). Choosing the branch of the logarithm with \(\log f(0)\in\mathbb R\), we observe that the radial limits of \(\log f\) exists for all \(\theta\in\mathbb T\setminus\{\theta_1,\dots, \theta_N\}\) and are given by 
\[
    \lim_{r\to 1^-}\log f(re^{i\theta})=-\sum_{j=1}^N\lambda_j\frac{1+e^{i(\theta-\theta_j)}}{1-e^{i(\theta-\theta_j)}}=-i\sum_{j=1}^N\lambda_j\cot\bigg(\frac{\theta-\theta_j}{2}\bigg).
\]
Defining 
\begin{equation}\label{eq: standard varphi}
    \varphi(\theta)=-\sum_{j=1}^N\lambda_j\cot\bigg(\frac{\theta-\theta_j}{2}\bigg),
\end{equation}
it follows that \(f(re^{i\theta})\to e^{i\varphi(\theta)}\) as \(r\to 1^-\) away from \(\theta_1,\dots,\theta_N\). The function \(e^{i\varphi(\theta)}\) therefore agrees a.e.\ with the non-tangential limit of \(f\), which exists a.e.\ in \(\mathbb T\) by Fatou's theorem \cite[Thm.\ 2.2]{Duren1970}. For convenience, we denote the a.e.\ defined non-tangential limit of \(f\) by \(f(e^{i\theta})\). 

Our main result is the following remarkable property of the function \(\varphi\). We just learnt that the first identity was proven by Poltoratski \cite{Poltoratski1996}, but we provide different proofs that we feel are of independent interest.

\begin{thm}\label{thm: hilbert identity}
Let \(\varphi\) be as in \eqref{eq: standard varphi}, and set \(t=\lambda_1+\cdots+\lambda_N\). Then for all \(y\in\mathbb R\) it holds that 
\[|\{\theta\in\mathbb T\;|\; \varphi(\theta)>y\}|=\pi-2\arctan(y/t),\] 
and
\[
    \sum_{\varphi(\theta)=y}\frac{1}{\varphi'(\theta)}=\frac{2t}{t^2+y^2}.    
\]
\end{thm}

Three proofs of this theorem are given, to illustrate its connection with various areas of complex and harmonic analysis. In a parallel work to this one \cite{MacDonald2026b}, the third author uses this formula in a fundamental way to prove a property of extremal functions in the Krzy\.z conjecture. Actually the current manuscript stems out of our attempts to prove such conjecture. The following consequence of Theorem \ref{thm: hilbert identity} is useful in computations. 

\begin{cor}\label{lem: change of variables formula}
Let \(F:\mathbb R\to\mathbb R\) be measurable. If one of the functions \(F^\pm=\pm\max\{\pm F,0\}\) satisfies \(F^{\pm}(x)/x^2\in L^1(\mathbb R)\), then 
\[
    \int_0^{2\pi}F(\varphi(\theta))\,d\theta=2t\int_{-\infty}^\infty\frac{F(x)}{x^2+t^2}\,dx.
\]
\end{cor}

Any function with negative real part in \(\mathbb D\) belongs to \(H^p(\mathbb D)\) for each \(p\in(0,1)\) \cite[Thm.\ 3.2]{Duren1970}, so \(\log f\in H^p(\mathbb D)\) for \(p\in(0,1)\). Using the preceding result, we give more detailed information. 

\begin{cor}\label{cor: H^p norm of log f}
If \(f\) is as in \eqref{ASIF} then \(\log f\not\in H^p(D)\) for \(p\geq 1\). If \(p\in(0,1)\), then 
\[
    \|\log f\|_{H^p(D)}=\frac{\lambda_1+\cdots+\lambda_N}{\cos(\pi p/2)^{1/p}}.
\]
\end{cor}

This calculation of norms allows us to give a completely new proof of the following well-known fact, which is usually derived from the properties of the Frostman shift and infinite Blaschke products.

\begin{cor}\label{cor: weighted sum is infinite}
Let \(f(z)=\sum_{j=0}^\infty a_jz^j\) be an atomic singular inner function. Then 
\[
    \sum_{j=1}^\infty j|a_j|^2=\infty.
\]
\end{cor}


By modifying our second proof of Theorem \ref{thm: hilbert identity}, we prove a result generalizing Corollary \ref{lem: change of variables formula}.

\begin{thm}\label{thm: upper half plane}
Let \(F\) be holomorphic in \(\mathbb H^+\), for \(x\in\mathbb R\) define \(\smash{F(x)\overset{a.e.}{=}\displaystyle\lim_{y\to 0^+}F(x+iy)}\), and assume that \(F(x)/x^2\in L^1(\mathbb R)\). If \(\varphi\) is as in \eqref{eq: standard varphi}, then 
\[
    \frac{1}{2\pi}\int_{0}^{2\pi}F(\varphi(\theta))d\theta=F(it),\qquad t=\lambda_1+\cdots+\lambda_N.
\]
\end{thm}

Finally, using Theorem \ref{thm: hilbert identity} and a simple rearrangement argument, we prove an interesting factorization of of the map \(\theta\mapsto f(e^{i\theta})\).

\begin{thm}\label{lem:boundaryfactorization}
Let \(f\) as in \eqref{ASIF}, and let \(t=\lambda_1+\cdots+\lambda_N\). Then there exists a map \(\sigma:\mathbb T\to\mathbb T\) with the following properties:
\begin{itemize}
    \item[(1)] If \(1\leq j\leq N\) then \(\sigma(\theta_j)=0\), and for all \(\theta\in\mathbb T\setminus\{\theta_1,\dots,\theta_N\}\) it holds that 
    \[
        f(e^{i\theta})\overset{a.e.}{=}\exp\bigg(-it\cot\bigg(\frac{\sigma(\theta)}{2}\bigg)\bigg),
    \]
    \item[(2)] The map \(\sigma\) is a continuous, monotone increasing, \(N\)-fold cover of \(\mathbb T\).
    \item[(3)] The map \(\sigma\) is measure-preserving. I.e.\,, if \(E\) is Lebesgue measurable then \(|\sigma^{-1}(E)|=|E|\). 
\end{itemize}
\end{thm}

\bigskip


\section{Preliminaries}

Given \(f:\mathbb T\to\mathbb C\), the Hilbert transform of \(f\) on the circle, denoted \(Hf\), is formally defined by the Cauchy principal value
\begin{equation}\label{def:HiblertTransform}
    Hf(x)=\frac{1}{2\pi}\mathrm{p.v.}\int_0^{2\pi}f(t)\cot\bigg(\frac{x-t}{2}\bigg)dt,\qquad x\in\mathbb T.
\end{equation}
Sometimes \(Hf\) is called the conjugate function to \(f\), since for \(1<p<\infty\) and real-valued \(f\in L^p(\mathbb T)\), the function \(f+iHf\) extends holomorphically to a function in \(H^p(D)\) by the Poisson integral. This and other properties of the operator \(H\) are well understood; see e.g. \cite[Chap.\ 2]{McGovern1980}. 

The following simple example will be useful in a moment.  If \(I=(a,b)\) is an interval in \(\mathbb T\), then the Hilbert transform of its indicator \(\chi_I\) is  given by
\begin{equation}\label{eq:Example HT of indicator}
    H\chi_I(x)=\frac{1}{2\pi}\mathrm{p.v.}\int_{a}^{b}\cot\bigg(\frac{x-t}{2}\bigg)dt=\frac{1}{\pi}\log\bigg|\frac{\sin(\frac{1}{2}(x-a))}{\sin(\frac{1}{2}(x-b))}\bigg|.
\end{equation}
In connection with this example, we have the following remarkable result of Stein \& Weiss \cite{SteinWeiss1959}.

\begin{thm}[Stein \& Weiss]\label{thm:SteinWeiss}
Let \(E\subseteq [0,2\pi)\) be a Lebesgue measurable set, and for \(y\in\mathbb R\) set \(\lambda_E^\pm(y)=|\{x\in \mathbb T\;|\; \pm H\chi_E(x)>y\}|\). Then
\begin{equation}\label{eq:SteinWeiss}
    e^{i\lambda^{\pm}( y)}=\frac{\sinh(\pi y)+i\sin(|E|/2)}{\sinh(\pi y)-i\sin(|E|/2)}.
\end{equation}
\end{thm}

\rmk Stein \& Weiss write \(\overline{\chi}_E\) instead of \(H\chi_E\). Their definition of the Hilbert transform differs ours one by a multiple of \(2\pi\), so we use the identity \(\overline{\chi}_E=2\pi H\chi_E\) when quoting \cite{SteinWeiss1959}.

The functions \(\lambda_E^\pm\), and the distribution function of \(H\chi_E\), are thus determined entirely by \(|E|\) independent of how the set is distributed. 

The following result functions as a quantitative version of the Borel-Cantelli lemma, and it will also be useful in what follows. We quote the following version from \cite[Thm.\ 1]{Yan2006}.

\begin{lem}[Kochen-Stone]\label{lem:KochenStone}
Let \(\{E_n\}_{n=1}^\infty\) be measurable subsets of \(\mathbb T\). If \(\sum_{n=1}^\infty|E_n|=\infty\), then
\[
    |\{x\in\mathbb T\;|\; x\in E_n\mathrm{\,\,infinitely\,\, often}\}|\geq \limsup_{n\to\infty}\frac{(\sum_{k=1}^n|E_n|)^2}{\sum_{k,l=1}^n|E_k\cap E_l|}.
\]
\end{lem}

\rmk This lemma is usually stated in terms of a probability measure, but the inequality is unchanged by rescaling the measure, so we can use it with Lebesgue measure on \(\mathbb T\).

\section{Proof of Theorem \ref{thm: hilbert identity} using the Hilbert transform}

Let \(\varphi\) be as in \eqref{eq: standard varphi}, and assume without loss of generality that \(0\leq \theta_1<\cdots<\theta_N<2\pi\). Momentarily fix \(\varepsilon>0\), and for each \(k\in\{1,\dots, N\}\) consider intervals \(I_k^{\varepsilon}\) in \(\mathbb T\) of length \(\lambda_k\varepsilon\) centred at each \(\theta_k\). Assume also that \(\varepsilon\) is small enough that these intervals are pairwise disjoint. Taking \(E^\varepsilon=\bigsqcup_{k=1}^nI_k^{\varepsilon}\) to be their disjoint union, the functions \(\chi_{E^{\varepsilon}}/\varepsilon\) approximate the density of the measure \(\sum_{k=1}^N\lambda_k\delta_{\theta_k}\) when \(\varepsilon\) is small. Furthermore, defining
\[
    \lambda^\pm_\varepsilon(y)=\bigg|\bigg\{x\in\mathbb T\;\bigg|\;\mp\frac{2\pi}{\varepsilon}H\chi_{E^{\varepsilon}}(x)>y\bigg\}\bigg|,\qquad y\in\mathbb R,
\]
we claim that for each \(y\in\mathbb R\), the limits \(\displaystyle\lim_{\varepsilon\to 0^+}\lambda^\pm_\varepsilon(y)\) exist and equal \(|\{x\in\mathbb T\;|\; \pm\varphi(x)>y\}|\).

To prove this, we first show that \(-2\pi H\chi_{E^{\varepsilon}}/\varepsilon\to \varphi\) pointwise on \(\mathbb T\setminus\{\theta_1,\dots,\theta_N\}\) as \(\varepsilon\to 0^+\). From \eqref{eq:Example HT of indicator}, it follows that
\[
    \lim_{\varepsilon\to 0^+}\frac{2\pi }{\varepsilon}H \chi_{E^{\varepsilon}}(x)=\lim_{\varepsilon\to 0^+}\frac{2\pi }{\varepsilon}\sum_{k=1}^NH \chi_{I_k^{\varepsilon}}(x)=\lim_{\varepsilon\to 0^+}\frac{2}{\varepsilon}\sum_{k=1}^N\log\bigg|\frac{\sin(\frac{1}{2}(x-\theta_k+\lambda_k\varepsilon/2))}{\sin(\frac{1}{2}(x-\theta_k-\lambda_k\varepsilon/2))}\bigg|.
\]
Using L'H\^opital's rule then gives for \(x\in \mathbb T\setminus\{\theta_1,\dots,\theta_N\}\)
\[
    \lim_{\varepsilon\to 0^+}\frac{2\pi }{\varepsilon}H \chi_{E^{(\varepsilon)}}(x)=\lim_{\varepsilon\to 0^+}\frac{1}{2}\sum_{k=1}^N\lambda_k\bigg(\cot\bigg(\frac{x-\theta_k+\lambda_k\varepsilon/2}{2}\bigg)+\cot\bigg(\frac{x-\theta_k-\lambda_k\varepsilon/2}{2}\bigg)\bigg)=-\varphi(x).
\]

Next, fix \(y\in\mathbb R\), set \(A^{\varepsilon}=\{x\in\mathbb T\;|\; -2\pi H\chi_{E^{\varepsilon}}(x)/\varepsilon>y\}\), and let \(A=\{x\in\mathbb T\;|\;\varphi(x)>y\}\). Our goal is to show that \(|A^\varepsilon|\to|A|\) as \(\varepsilon\to 0^+\). Equivalently, \(|A^\varepsilon\Delta A|\to 0\) as \(\varepsilon\to 0^+\), where \(\Delta\) denotes the symmetric difference. Assume toward a contradiction that there exists a sequence \(\varepsilon_n\to 0\) such that \(|A^{\varepsilon_n}\Delta A|\not\to 0\). Passing so a subsequence if necessary, we can assume that there exists \(m>0\) such that \(|A^{\varepsilon_n}\Delta A|\to m\). Defining \(B_n=A^{\varepsilon_n}\setminus A\) and \(C_n=A\setminus A^{\varepsilon_n}\), we point out that \(A^{\varepsilon_n}\Delta A=B_n\sqcup C_n\) so either \(|B_n|>m/2\) for \(n\) sufficiently large, or \(|C_n|>m/2\) for large \(n\). 

In the first case, we pass to a subsequence to assume that \(|B_n|>m/2\) for all \(n\in\mathbb N\), so that
\[
    \sum_{n=0}^\infty |B_n|=\infty.
\]
It follows from the Kochen-Stone lemma (Lemma \ref{lem:KochenStone} above) that
\begin{align*}
    |\{x\in\mathbb T\;|\; x\in B_n\textrm{ for infinitely many }n\}|&\geq \limsup_{n\to\infty}\frac{(\sum_{k=1}^n|B_k|)^2}{\sum_{k,l=1}^n|B_k\cap B_l|}\\
    &\geq \limsup_{n\to\infty}\frac{(\sum_{k=1}^n|B_k|)^2}{2\pi n^2}\geq \frac{m^2}{8\pi }.
\end{align*}
For the latter estimates we used that \(|B_k\cap B_l|\leq 2\pi \) and that \(|B_k|\geq m/2\). If \(x\in B_n\) for infinitely many \(n\) and \(x\not\in \{\theta_1,\dots,\theta_n\}\) then there exists a subsequence \(\varepsilon_{n_k}\to 0\) such that for all \(k\),
\[
    \varphi(x)\leq y<-\frac{2\pi}{\varepsilon_{n_k}}H\chi_{E^{\varepsilon_{n_k}}}(x).
\]
Sending \(k\to\infty\), it follows that \(\varphi(x)=y\). As \(\varphi\) is monotone on each interval \((\theta_{j},\theta_{j+1})\) this equation has at most \(N\) solutions, so 
\[
    \#\{x\in\mathbb T\;|\; x\in B_n\textrm{ for infinitely many }n\}\leq \#\{\theta_1,\dots,\theta_N\}+\#\{x\in\mathbb T\;|\;\varphi(x)=y\}\leq 2N. 
\]
In particular the set on the left is finite, a contradiction since it has positive Lebesgue measure.

If instead \(|C_n|>m/2\) for large \(n\), then once again we can assume it holds for all \(n\in\mathbb N\) by passing to a subsequence. It follows that \(|\{x\in\mathbb T\;|\; x\in C_n\textrm{ for infinitely many }n\}|>0\) by an identical argument. If there exists \(x\in\mathbb T\) belonging to this set and \(x\not\in\{\theta_1,\dots,\theta_N\}\) then there is a subsequence \(\varepsilon_{n_k}\to 0\) such that for all \(k\),
\[
    \varphi(x)>y\geq -\frac{2\pi}{\varepsilon_{n_k}}H\chi_{E^{\varepsilon_{n_k}}}(x).
\]
Sending \(k\to\infty\) yields \(\varphi(x)>\varphi(x)\), a contradiction. It follows that
\[
    \{x\in\mathbb T\;|\; x\in C_n\textrm{ for infinitely many }n\}\subseteq \{\theta_1,\dots,\theta_N\},
\]
and again we have found a finite set with positive Lebesgue measure, a contradiction. 

It follows now that \(|A^\varepsilon\Delta A|\to 0\) as \(\varepsilon\to 0^+\), meaning that \(\lambda^+_\varepsilon(y)\to  |\{x\in\mathbb T\;|\; \varphi(x)>y\}|\) as \(\varepsilon\to 0^+\). The claim for \(\lambda^-_\varepsilon\) follows from recognizing that \(\lambda_{\varepsilon}^-(y)=2\pi-\lambda^+_\varepsilon(-y)\) for all \(\varepsilon>0\) and \(y\in\mathbb R\), since \(H\chi_{E^\varepsilon}(x)=y\) has only finitely many solutions for each \(y\in\mathbb R\). Likewise,
\[
    |\{x\in\mathbb T\;|\; \varphi(x)<-y\}|=2\pi -|\{x\in\mathbb T\;|\; \varphi(x)>-y\}|,
\]
and the convergence \(\lambda_\varepsilon^-(y)\to |\{x\in\mathbb T\;|\; -\varphi(x)>y\}|\) is immediate.

Equipped with this convergence, we use Theorem \ref{thm:SteinWeiss} to compute \(\lambda_\varepsilon(y)\) for \(\varepsilon>0\) before taking a limit as \(\varepsilon\to 0^+\). Write \(t=\lambda_1+\cdots+\lambda_N\) so that \(|E^{\varepsilon}|=\sum_{k=1}^N|I_k^{\varepsilon}|=\sum_{k=1}^N\lambda_k\varepsilon=\varepsilon t\). Equation \eqref{eq:SteinWeiss} shows that for any \(y\in\mathbb R\),
\[
    e^{i\lambda^-_\varepsilon(y)}=e^{i|\{H_\mathbb T\chi_{E^{\varepsilon}}>\varepsilon y\}|}=\frac{\sinh(\varepsilon y/2)+i\sin(\varepsilon t/2)}{\sinh(\varepsilon y/2)-i\sin(\varepsilon t/2)}.
\]
Sending \(\varepsilon\to 0^+\) it follows from the convergence of the distribution functions \(\lambda^{\pm}_\varepsilon\) established above that if \(y\neq 0\) then 
\[
    e^{i|\{x\in\mathbb T\;|\; -\varphi(x)>y\}|}=\lim_{\varepsilon\to 0^+}e^{i\lambda^-_\varepsilon(y)}=\lim_{\varepsilon\to 0^+}\frac{1+i\sin(\varepsilon t/2)/\sinh(\varepsilon y/2)}{1-i\sin(\varepsilon t/2)/\sinh(\varepsilon y/2)}=\frac{y+it}{y-it}.
\]
The constant on the right is unimodular, so taking a logarithm with argument in \([0,2\pi)\) yields
\[
    |\{x\in\mathbb T\;|\;\varphi(x)<-y\}|=|\{x\in\mathbb T\;|\;-\varphi(x)>y\}|=\mathrm{arg}\bigg(\frac{y+it}{y-it}\bigg)=\pi-2\arctan\bigg(\frac{y}{t}\bigg).
\]
As \(|\{x\in\mathbb T\;|\;\varphi(x)=-y\}|=0\) since \(\varphi\) is piecewise monotone, we see that for any \(y\neq 0\),
\[
    |\{x\in\mathbb T\;|\;\varphi(x)>y\}|=2\pi-|\{x\in\mathbb T\;|\;\varphi(x)<y\}|=\pi-2\arctan\bigg(\frac{y}{t}\bigg).
\]
Finally, if \(y=0\) then \(e^{i|\{x\in\mathbb T\;|\; \varphi(x)<0\}|}=-1\), so \(|\{x\in\mathbb T\;|\; \varphi(x)<0\}|=\pi\) and \(|\{x\in\mathbb T\;|\; \varphi(x)>0\}|=\pi\), meaning that the formula above also holds when \(y=0\). 

Next, for \(1\leq j\leq N\) define \(\varphi_j:=\varphi|_{(\theta_j,\theta_{j+1})}\). Each of these functions increases monotonically from \(-\infty\) to \(\infty\), so they invertible with monotone inverses \(\varphi_j^{-1}:\mathbb R\to(\theta_j,\theta_{j+1})\). Indeed, for each \(y\in\mathbb R\) it holds that \(\varphi_j^{-1}(y)=\theta_{j+1}-\theta_j-|\{\theta\in\mathbb T\;|\; \varphi(\theta)>y\}\cap(\theta_j,\theta_{j+1})|\). Summing over \(j\) gives
\[
    \sum_{j=1}^N\varphi_j^{-1}(y)=2\pi-|\{\theta\in\mathbb T\;|\; \varphi(\theta)>y\}|=\pi+2\arctan\left(\frac{y}{t}\right),
\]
and differentiating in \(y\) in turn yields
\begin{equation}\label{eq: Jacobian formula for varphi}
    \sum_{j=1}^N\frac{1}{\varphi'(\varphi_j^{-1}(y))}=\frac{2t}{y^{2}+t^{2}}.    
\end{equation}
The numbers \(\{\varphi_1^{-1}(y),\dots,\varphi_N^{-1}(y)\}\) are the solutions to the equation \(\varphi(\theta)=y\), so we are done. \qed

\section{Alternate proof using the mean value property of harmonic functions}

Let \(\theta_1,\dots,\theta_N\) be as in \eqref{eq: standard varphi}, write \(t=\lambda_1+\cdots+\lambda_N\), and define a function \(\Phi\in \mathrm{Hol}(\mathbb D)\) by
\begin{equation}\label{eq: Phi}
    \Phi(z)=-i\sum_{j=1}^N\lambda_j\frac{z+e^{i\theta_j}}{z-e^{i\theta_j}}.
\end{equation}
On the circle, we note that \(\Phi(e^{i\theta})=\varphi(\theta)\) away from \(\{\theta_1,\dots,\theta_N\}\). Furthermore, \(\im \Phi=0\) a.e. on \(\partial \mathbb D\) and \(\im\Phi (0)=t>0\), so \(\im\Phi>0\) in \(\mathbb D\) by the open mapping theorem. For fixed \(y\in\mathbb R\) we have \(\im \{\Phi-y\}=\im \Phi>0\) in \(\mathbb D\), so the function \(\Psi_y(z)=\log(\Phi(z)-y)\) (defined using the principal logarithm) is holomorphic in \(\mathbb D\). Thus, \(\im \Psi_y\) is harmonic and takes values in \((-\pi,\pi]\). 

By the mean value property of harmonic functions, we observe that for every \(r\in(0,1)\),
\[
    \frac{1}{2\pi}\int_0^{2\pi}\im\Psi_y(re^{i\theta})d\theta=\im \Psi_y(0)=\im \log(it-y)=\arg(it-y)=\frac{\pi}{2}+\arctan\bigg(\frac{y}{t}\bigg).
\]
For \(\theta\not\in\{\theta_1,\dots,\theta_N\}\) also note that \(\Phi(re^{i\theta})\to \varphi(\theta)\in\mathbb R\) as \(r\to 1^-\), so \(\Phi(re^{i\theta})-y\) approaches the real axis as \(r\to 1^-\) for a.e. \(\theta\in\mathbb T\). Depending on the sign of the limit, we thus find that
\[
    \lim_{r\to 1^-}\im\Psi_y(re^{i\theta})=\lim_{r\to 1^-}\arg(\Phi(re^{i\theta})-y)\overset{a.e.}{=}\begin{cases}
        0 & \varphi(\theta)\geq y\\
        \pi & \varphi(\theta)< y. 
    \end{cases}
\]
It follows from the dominated convergence theorem that
\[
    \pi+2\arctan\bigg(\frac{y}{t}\bigg)=\frac{1}{\pi}\lim_{r\to 1^-}\int_0^{2\pi}\im\Psi_y(re^{i\theta})d\theta=\frac{1}{\pi}\int_0^{2\pi}\lim_{r\to 1}\im\Psi_y(re^{i\theta})d\theta=\int_{0}^{2\pi}\chi_{\{\varphi(\theta)<y\}}d\theta.
\]
The right-hand side is exactly \(|\{\theta\in\mathbb T\;|\; \varphi(\theta)<y\}|\). As \(\varphi\) is strictly increasing, for each \(y\in \mathbb R\) the set \(\{\theta\in\mathbb T\;|\;\varphi(\theta)=y\}\) is finite, hence it has zero measure. Therefore, 
\[
    |\{\theta\in\mathbb T\;|\; \varphi(\theta)>y\}|=2\pi -|\{\theta\in\mathbb T\;|\; \varphi(\theta)<y\}|=2\pi -\bigg(\pi+2\arctan\bigg(\frac{y}{t}\bigg)\bigg)=\pi-2\arctan\bigg(\frac{y}{t}\bigg).
\]
The remainder of the proof proceeds exactly as above.

\section{Alternate proof using Blaschke products}

It has long been known that if \(h\) is a finite Blaschke product and \(h(0)=0\), then \(h|_{\partial \mathbb D}:\partial\mathbb D\to \partial \mathbb D\) is measure preserving with respect to Lebesgue measure: see, for instance the more general result of Nordgren \cite[Cor.\ on p.\ 442]{Nordgren_1968}. In \cite[\S6]{MacDonald2026}, it is shown that if \(f\) is an atomic singular inner function as in \eqref{ASIF} with \(f(0)>0\), and if \(t=\lambda_1+\cdots+\lambda_n\), then there exists a degree \(N\) Blaschke product \(h\) such that \(h(0)=0\) and 
\[
    f(z)=\exp\bigg(t\frac{h(z)-1}{h(z)+1}\bigg),\qquad z\in D. 
\]
Taking radial limits, it follows at once that away from \(\theta_1,\dots,\theta_N\),
\[
    i\varphi(\theta)=\log f(e^{i\theta})=t\frac{h(e^{i\theta})-1}{h(e^{i\theta})+1}=\frac{2it\im\{h(e^{i\theta})\}}{|h(e^{i\theta})+1|^2}.
\]
Writing \(h(e^{i\theta})=e^{i\eta(\theta)}\) for \(\eta\) taking values in \(\mathbb R\), this gives 
\[
    \varphi(\theta)=\frac{2t\im\{e^{i\eta(\theta)}\}}{|e^{i\eta(\theta)}+1|^2}=\frac{2t\sin(\eta(\theta))}{2+2\cos(\eta(\theta))}=t\tan\bigg(\frac{\eta(\theta)}{2}\bigg). 
\]
As \(\tan(\theta/2)\) is strictly increasing on \((-\pi,\pi)\) and the above holds except at finitely many points, for any \(y\in\mathbb R\) we have
\[
    |\{\theta\in\mathbb T\;|\; \varphi(\theta)>y\}|=\bigg|\bigg\{\theta\in\mathbb T\;\bigg|\; t\tan\bigg(\frac{\eta(\theta)}{2}\bigg)>y\bigg\}\bigg|=\bigg|\bigg\{\theta\in\mathbb T\;\bigg|\; \eta(\theta)>2\arctan\bigg(\frac{y}{t}\bigg)\bigg\}\bigg|. 
\]
Finally, we note that as \(h\) is measure-preserving on the circle, the function \(\eta\) is uniformly distributed over \(\mathbb T\), meaning that
\[
    |\{\theta\in\mathbb T\;|\; \varphi(\theta)>y\}|=\bigg|\bigg\{\theta\in\mathbb T\;\bigg|\; \theta>2\arctan\bigg(\frac{y}{t}\bigg)\bigg\}\bigg|=\pi-2\arctan\bigg(\frac{y}{t}\bigg).
\]
Again, the proof from here proceeds identically.

\section{Proofs of Corollaries}

\begin{proof}[Proof of Corollary \ref{lem: change of variables formula}]
On each interval \((\theta_j,\theta_{j+1})\) we perform the simple change of variables \(x=\varphi_j^{-1}(u)\), which takes \((\theta_j,\theta_{j+1})\) to \(\mathbb R\). As \(dx=\varphi_j'(\varphi_j^{-1}(u))^{-1}du\) we have
\[
    \int_0^{2\pi}F(\varphi(x))\,dx=\sum_{j=1}^N\int_{\theta_j}^{\theta_{j+1}}F(\varphi_j(x))\,dx=\sum_{j=1}^N\int_{-\infty}^{\infty}\frac{F(u)}{\varphi_j'(\varphi_j^{-1}(u))}\,du.
\]
The sum above is well-defined thanks to our condition on \(F^{\pm}\) (if any integral on the right takes the value \(\pm\infty\), then none of the other terms can be \(\mp\infty\)). Thus, by \eqref{eq: Jacobian formula for varphi},
\[
    \int_0^{2\pi}F(\varphi(x))\,dx=\int_{-\infty}^{\infty}F(u)\bigg(\sum_{j=1}^N\frac{1}{\varphi_j'(\varphi_j^{-1}(u))}\bigg)\,du=2t\int_{-\infty}^{\infty}\frac{F(u)}{u^2+t^2}\,du.
\]
This is the claimed formula.
\end{proof}

\begin{proof}[Proof of Corollary \ref{cor: H^p norm of log f}]
Recall (see e.g. \cite[Thm. 3.8]{Katznelson2004}) that if \(p>0\) and if \(g\in\mathrm{Hol}(D)\) has radial limit \(\tilde{g}(\theta)\overset{a.e.}{=}g(e^{i\theta})\), then \(\|g\|_{H^p(D)}^p=\|\tilde{g}\|_{L^p(\mathbb T)}^p/2\pi\). Taking \(g=\log f\), it follows that \(\tilde{g}=i\varphi\) and
\[
    \|\log f\|_{H^p(D)}^p=\frac{1}{2\pi}\|\varphi\|_{L^p(\mathbb T)}^p=\frac{1}{2\pi}\int_0^{2\pi} |\varphi(x)|^pdx.
\]
Taking \(F(x)=|x|^p\), it is clear that \(F\geq 0\) and \(F^-\equiv 0\in L^1(\mathbb R)\), so by Corollary \ref{lem: change of variables formula},
\[
    \|\log f\|_{H^p(D)}^p=\frac{t}{\pi}\int_{-\infty}^\infty\frac{|u|^p}{u^2+t^2}du=\frac{2t}{\pi}\int_{0}^\infty\frac{u^p}{u^2+t^2}du=\frac{2t^{p}}{\pi}\int_{0}^\infty\frac{u^p}{u^2+1}du.
\]
The integrand on the right is comparable to \( u^{p-1}\) when \(u\) is large, so when \(p\geq 1\) it diverges and \(\log f\not\in H^p(D)\). If \(0<p<1\) then the integral converges, and the claimed value follows from a standard contour integration.
\end{proof}

\begin{proof}[Proof of Corollary \ref{cor: weighted sum is infinite}]
Writing \(f=e^g\), we have \(zf'(z)\overline{f(z)}=zg'(z)|f(z)|^2\). Thus by Parseval's identity, for each \(r\in(0,1)\) it holds that
\[
    \sum_{j=1}^\infty j|a_j|^2\geq \sum_{j=1}^\infty j|a_j|^2r^{2j}=\frac{1}{2\pi}\int_0^{2\pi}re^{i\theta}f'(re^{i\theta})\overline{f(re^{i\theta})}d\theta=\frac{1}{2\pi}\int_0^{2\pi}re^{i\theta}g'(re^{i\theta})|f(re^{i\theta})|^2d\theta.
\]
Sending \(r\to 1^-\), Fatou's lemma shows that
\[
    \sum_{j=1}^\infty j|a_j|^2\geq \liminf_{r\to 1}\frac{1}{2\pi}\int_0^{2\pi}re^{i\theta}g'(re^{i\theta})|f(re^{i\theta})|^2d\theta\geq \frac{1}{2\pi}\int_0^{2\pi}\liminf_{r\to 1}re^{i\theta}g'(re^{i\theta})|f(re^{i\theta})|^2d\theta.
\]
The right-hand side is exactly \(\int_0^{2\pi}\varphi'(\theta)d\theta/2\pi\). Observing that \(\varphi'\geq |\varphi|\) pointwise on \(\mathbb T\setminus\{\theta_1,\dots,\theta_N\}\), it follows that \(\sum_{j=1}^\infty j|a_j|^2\geq \|\varphi\|_{L^1(\mathbb T)}= \|\log f\|_{H^1(D)}=\infty\) as desired.
\end{proof}

\section{Proof of Theorem \ref{thm: upper half plane}}

Let \(\Phi\) be as in \eqref{eq: Phi}, and observe that for any \(z\in D\) it holds that
\[
    \im \Phi(z)=\im\bigg\{-i\sum_{j=1}^N\lambda_j\frac{z+e^{i\theta_j}}{z-e^{i\theta_j}}\bigg\}=-\sum_{j=1}^N\lambda_j\re\bigg\{\frac{z+e^{i\theta_j}}{z-e^{i\theta_j}}\bigg\}=\sum_{j=1}^N\lambda_j\frac{1-|z|^2}{|z-e^{i\theta_j}|^2}>0.
\]
Thus, \(G:=F\circ\Phi\) is well-defined and holomorphic in \(\mathbb D\). Cauchy's integral formula  also gives
\[
    \frac{1}{2\pi}\int_0^{2\pi }G(re^{i\theta})d\theta=G(0)=F(\Phi(0))=F(it),\qquad r\in(0,1)
\]
For \(r<1\) and \(\theta\in\mathbb T\) define \(\gamma_r(\theta)=G(re^{i\theta})\) and \(\gamma(\theta)=G(e^{i\theta})\). Then \(\gamma_r\to \gamma\) pointwise a.e. in \(\mathbb T\) as \(r\to 1^-\), and it also follows from the definition of \(G\) that \(\gamma=F\circ\varphi\) for a.e. \(\theta\in\mathbb T\). Our integrability condition on \(F\) and Corollary \ref{lem: change of variables formula} ensure that \(\gamma=F\circ\varphi\in L^1(\mathbb R)\), and by Hardy's convexity theorem (see e.g.\ \cite[Thm. 1.5]{Duren1970}), the map \(r\to \int_0^{2\pi}|G(re^{i\theta})|\,d\theta/2\pi=\|\gamma_r\|_{L^1(\mathbb T)}\) is non-decreasing. Therefore, \(\|\gamma_r\|_{L^1(\mathbb T)}\to \|\gamma\|_{L^1(\mathbb T)}\) as \(r\to 1^-\). 

Using the aforementioned properties of the functions \(\gamma_r\), a standard argument will yield the desired formula. Consider the functions \(\eta_r=|\gamma_r|+|\gamma|-|\gamma-\gamma_r|\), which satisfy \(\eta_r\geq 0\) and \(\eta_r\to 2|g|\) a.e. as \(r\to 1^-\). Thus by Fatou's lemma, 
\[
    \|\gamma\|_{L^1(\mathbb T)}=\frac{1}{4\pi}\int_0^{2\pi}\lim_{r\to 1^-}\eta_r(\theta)\,d\theta\leq  \frac{1}{4\pi}\liminf_{r\to 1^-}\int_0^{2\pi}\bigg(|\gamma_r(\theta)|+|\gamma(\theta)|-|\gamma(\theta)-\gamma_r(\theta)|\bigg)d\theta.
\]
Using convergence of the \(L^1(\mathbb T)\) norms, we find that
\[
    \|\gamma\|_{L^1(\mathbb T)}\leq \|\gamma\|_{L^1(\mathbb T)} -\frac{1}{4\pi}\limsup_{r\to 1^-}\int_0^{2\pi}|\gamma(\theta)-\gamma_r(\theta)|d\theta,
\]
and since \(\|\gamma\|_{L^1(\mathbb T)}<\infty\) we deduce that the limsup on the right is identically zero. It follows that \(\|\gamma_r-\gamma\|_{L^1(\mathbb T)}\to 0\) as \(r\to 1^-\), meaning that
\[
    \limsup_{r\to 1^-}\bigg|\int_0^{2\pi }(G(re^{i\theta})-G(e^{i\theta}))d\theta\bigg|\leq \limsup_{r\to 1^-}\int_{0}^{2\pi }|G(re^{i\theta})-G(e^{i\theta})|d\theta=2\pi\lim_{r\to 1^-}\|\gamma_r-\gamma\|_{L^1(\mathbb T)},
\]
and the limit on the right is zero. The claimed formula now follows. \qed

\section{Proof of Theorem \ref{lem:boundaryfactorization}}

For convenience, set \(\lambda(y)=|\{x\in\mathbb T\;|\;\varphi(x)>y \}|\). We define the monotone non-decreasing rearrangement of \(\varphi\), denoted \(\varphi^*:\mathbb T\to\mathbb R\), by \(\varphi^*(\theta)=\inf\{y\in\mathbb R\;|\; \lambda(y)\leq 2\pi -\theta\}\). Using Theorem \ref{thm: hilbert identity} we can compute this explicitly,
\[
    \varphi^*(\theta)=\inf\bigg\{y\in\mathbb R\;|\; \arctan\left(\frac{y}{t}\right)\geq \frac{\theta-\pi}{2}\bigg\}=\inf\bigg\{y\in\mathbb R\;|\; y\geq t\tan\bigg(\frac{\theta-\pi}{2}\bigg)\bigg\}=-t\cot\bigg(\frac{\theta}{2}\bigg).
\]
By construction, \(|\{x\in\mathbb T\;|\;\varphi^*(x)>y\}|=\lambda(y)\). Next we define \(\sigma:\mathbb T\to\mathbb T\) by \(\sigma(\theta)=2\pi-\lambda(\varphi(\theta))\). Since \(\lambda\) is continuous and monotone, for any \(\theta\in\mathbb T\) it holds that
\[
    \varphi^*(\sigma(\theta))=\inf\{y\in\mathbb R\;|\; \lambda(y)\leq 2\pi -\sigma(\theta)\}=\inf\{y\in\mathbb R\;|\; \lambda(y)\leq \lambda(\varphi(\theta))\}=\varphi(\theta).
\]
The representation of \(f\) in \textit{(1)} follows immediately.

It is clear that \(\sigma\) is smooth away from \(\{\theta_1,\dots,\theta_N\}\), and monotonicity follows from the bound
\[
    \sigma(\theta)=2\pi-\lambda(\varphi(\theta))=\pi+2\arctan\left(\frac{\varphi(\theta)}{t}\right)\qquad \rightarrow\qquad \sigma'(\theta)=\frac{2t\varphi'(\theta)}{t^2+\varphi(\theta)^2}\geq1.
\]
Moreover, by defining \(\sigma(\theta_j)=0\) for \(1\leq j\leq N\) extends \(\sigma\) to a continuous function on \(\mathbb T\) since
\[
    \lim_{\theta\to \theta_j}\sigma(\theta)=2\pi-\lim_{\theta\to \theta_j}\lambda(\varphi(\theta))=2\pi-\lim_{y\to \pm\infty}\lambda(y)\equiv 0\mod 2\pi.
\]

It remains to verify that \(\sigma\) is measure-preserving. It suffices to show that \(|\sigma^{-1}([0,\alpha])|=\alpha\) for any \(\alpha\in[0,2\pi]\), since the intervals \([0,\alpha]\) generate the Borel \(\sigma\)-algebra on \(\mathbb T\) and since every Lebesgue measurable set can be written as the union of a Borel set and a subset of a null Borel set. By construction, \(\sigma\geq 0\), so we can write \(\sigma^{-1}([0,\alpha])=\{\theta\in \mathbb T\;|\; \sigma(\theta)\leq \alpha\}\). The inequality \(\sigma(\theta)\leq \alpha\) can be written
\[
    \arctan\left(\frac{\varphi(\theta)}{t}\right)\leq \frac{\alpha-\pi}{2},
\]
and since the tangent function is strictly increasing on \((-\pi/2,\pi/2)\) and \(t>0\), the above holds if and only if 
\[
    \varphi(\theta)\leq t\tan\bigg(\frac{\alpha-\pi}{2}\bigg)=-t\cot\bigg(\frac{\alpha}{2}\bigg).
\]
It follows that
\[
    \sigma^{-1}([0,\alpha])=\bigg\{\theta\in[0,2\pi)\;\bigg|\;\varphi(\theta)\leq -t\cot\bigg(\frac{\alpha}{2}\bigg)\bigg\}=2\pi\setminus\bigg\{\theta\in[0,2\pi)\;\bigg|\;\varphi(\theta)> -t\cot\bigg(\frac{\alpha}{2}\bigg)\bigg\},
\]
and \(|\sigma^{-1}([0,\alpha])|=\sigma(-t\cot(\frac{\alpha}{2}))=\pi+2\arctan(-\cot(\frac{\alpha}{2}))=\alpha\), as desired. \qed

\subsection*{Acknowledgments} 

The authors would like to thank Almut Burchard, Fedor Nazarov, Kirill Serkh, and Stefan Steinerberger for their interest in this research and for helpful conversations, comments, and suggestions on the subject. We are grateful for the ideas and perspectives they shared.

Efraimidis, Llinares, Uriarte-Tuero and Vukoti\'c are partially supported by PID2025-168349NB-I00 funded by MICIU /AEI /10.13039/501100011033 / FEDER, UE (ERDF funds, EU). 

Eifraimidis and Liinares are  also supported by CIGE/2026/195 from Conselleria d’Edu\-caci\'o, Cultura i Universitats, Generalitat Valenciana, Spain. 

Llinares is additionally supported by PID2024-160326NA-I00, also  from MICIU. 

MacDonald is supported by an NSERC CGS scholarship and graduate funding from the University of Toronto. 

Uriarte-Tuero is also partially supported by an NSERC Discovery Grant, Canada. 

\subsection*{AI Disclosure} 
AI was only used in the proof-reading of the manuscript, which provided further references we were not aware of, but all the mathematical proofs and writing were done without AI.

\bibliographystyle{siam}
\bibliography{references}
\address

\end{document}